\input amstex
\documentstyle{amsppt}
\magnification 1200
\vcorrection{-9mm}
\input epsf
\NoBlackBoxes

{\hfill\it To Vladimir Lin in occasion of his 85-th birthday}
\bigskip

\topmatter
\title Homomorphisms of commutator subgroups of braid groups
       with small number of strings
\endtitle

\author S.~Yu.~Orevkov
\endauthor

\address IMT, Univ.~Paul Sabatier, Toulouse, France;
         Steklov Math.~Inst., Moscow, Russia 
\endaddress

\abstract
For any $n$, we describe all endomorphisms of
the braid group $B_n$ and of its commutator subgroup $B'_n$,
as well as all homomorphisms $B'_n\to B_n$.
These results are new only for small $n$ because
endomorphisms of $B_n$ are already described by Castel for $n\ge 6$,
and homomorphisms $B'_n\to B_n$ and endomorphisms of $B'_n$
are already described by Kordek and Margalit for $n\ge 7$.
We use very different approaches for $n=4$ and for $n\ge 5$.
\endabstract

\endtopmatter
\rightheadtext{Homomorphisms of commutator subgroups of braid groups}

\def\B{\bold B}
\def\P{\bold P}
\def\J{\bold J}
\def\A{\bold A}
\def\K{\bold K}
\def\bS{\bold S}

\def\C{\Bbb C}
\def\Z{\Bbb Z}

\def\ab{{\frak{ab}}}
\def\id{\operatorname{id}}
\def\im{\operatorname{im}}
\def\lk{\operatorname{lk}}
\def\sh{\operatorname{sh}}
\def\Sym{\operatorname{Sym}}
\def\Aut{\operatorname{Aut}}
\def\Out{\operatorname{Out}}
\def\Inn{\operatorname{Inn}}
\def\GL{\operatorname{GL}}
\def\SL{\operatorname{SL}}

\def\sectBn    {2}
\def\sectBfour {3}

\def\thBBn     {1.1}
\def\corBBnOne {1.2}
\def\corBBnTwo {1.3}
\def\thBBfour  {1.4}
\def\corBBfour {1.5}
\def\thBn      {1.6}
\def\thBfour   {1.7}
\def\propBthree  {1.8}
\def\propBthreeBfour {1.9}
\def\remTransv {1.10}

\def\lemPab   {\sectBn.1}
\def\lemJab   {\sectBn.2}
\def\lemMu    {\sectBn.3}
\def\propKM   {\sectBn.4}
\def\lemKMsix {\sectBn.5}
\def\lemKM    {\sectBn.6}
\def\lemXYX   {\sectBn.7}
\def\lemOne   {\sectBn.8}
\def\lemThree {\sectBn.9}
\def\lemSeven {\sectBn.10}
\def\lemEight {\sectBn.11}
 \let\remGarside=\remG

\def\lemInj   {\sectBfour.1}
\def\lemCommCG{\sectBfour.2}
\def\lemG     {\sectBfour.3}
\def\lemCD    {\sectBfour.4}
\def\lemCDii  {\sectBfour.5}
\def\lemK     {\sectBfour.6}
\def\lemF     {\sectBfour.7}
\def\lemXF    {\sectBfour.8}
\def\lemWF    {\sectBfour.9}
\def\lemUT    {\sectBfour.10}
\def\lemR     {\sectBfour.11}

\def\eqDefR  {1}

\def\eqKlein {2}
\def\eqKM    {3}
\def\eqLK    {4}

\def\eqGL    {5}
\def\eqGLsig {6}
\def\eqDTU   {7}
\def\eqCGC   {8}
\def\eqCWC   {9}
\def\eqPXP   {10}
\def\eqPGP   {11}
\def\eqId    {12}
\def\eqLemUT {13}

\def\figDTU  {1}
\def\figCGCC {2}
\def\figCG   {3}
\def\figCRS  {4}

\def\refArtin  {1}
\def\refBM     {2} 
\def\refBLM    {3} 
\def\refCastel {4}
\def\refDG     {5} 
\def\refFM     {6} 
\def\refGM     {7} 
\def\refGMW    {8} 
\def\refGL     {9} 
\def\refI      {10} 
\def\refKM     {11} 
\def\refLin    {12}
\def\refLinFAA {13}
\def\refLinTalk{14}
\def\refMKS    {15} 
\def\refJAlg   {16}
\def\refAFST   {17}
\document

\head Introduction
\endhead

Let $\B_n$ be the braid group with $n$ strings.
It is generated by $\sigma_1,\dots,\sigma_{n-1}$ (called {\it standard}
or {\it Artin} generators) subject to the relations
$$
   \text{$\sigma_i\sigma_j=\sigma_j\sigma_i$ for $|i-j|>1$;}\qquad
   \text{$\sigma_i\sigma_j\sigma_i=\sigma_j\sigma_i\sigma_j$ for
             $|i-j|=1$}.
$$
Let $\B'_n$ be the commutator subgroup of $\B_n$.

In this paper we describe all endomorphisms of $\B_n$ and $\B'_n$ and
homomorphisms $\B'_n\to\B_n$ for any $n$.
These results are new only for small $n$ because
endomorphisms of $\B_n$ are described by Castel in [\refCastel] for $n\ge 6$,
and homomorphisms $\B'_n\to\B_n$ and endomorphisms of $\B'_n$
are described by Kordek and Margalit in [\refKM] for $n\ge 7$.

The automorphisms of $\B_n$ and $\B'_n$ have been already known for any $n$:
Dyer and Grossman [\refDG] proved
that the only non-trivial element of $\Out(\B_n)$ corresponds to the automorphism
$\Lambda$ defined by $\sigma_i\mapsto\sigma_i^{-1}$ for any $i=1,\dots,n-1$, and
in [\refAFST] we proved that the restriction map 
$\Aut(\B_n)\to\Aut(\B'_n)$ is an isomorphism for $n\ge 4$ ($\B'_3$ is
a free group of rank $2$, thus its automorphisms are known as well; see e.g.~[\refMKS]).

The problem to study homomorphisms between braid groups and, especially,
between their commutator subgroups was posed by Vladimir Lin
[\refLin--\refLinTalk] because he found its applications to
the problem of superpositions of algebraic functions (the initial motivation for
Hilbert's 13th problem), see [\refLinFAA] and references therein.

Let us formulate the main results. We start with those about
homomorphisms of $\B'_n$ to $\B_n$ and to itself.

\proclaim{ Theorem \thBBn} {\rm(proven for $n\ge 7$ in [\refKM]).} 
Let $n\ge 5$. Then every non-trivial homomorphism $\B'_n\to\B_n$
extends to an automorphism of $\B_n$.
\endproclaim

We proof this theorem in \S\sectBn.
Since $\B''_n=\B'_n$ and $\Aut(\B_n)=\Aut(\B'_n)$ for $n\ge 5$,
the following two corollaries are, in fact, equivalent versions of
Theorem~\thBBn.

\proclaim{ Corollary \corBBnOne }
If $n\ge 5$, then any non-trivial endomorphism of $\B'_n$
is bijective.
\endproclaim

\proclaim{ Corollary \corBBnTwo }
If $n\ge 5$, then any non-trivial homomorphism $\B'_n\to\B_n$
is an automorphism of $\B'_n$ composed with the
inclusion map.
\endproclaim


Let $R$ be the homomorphism
$$
    R:\B_4\to\B_3, \qquad
    \sigma_1,\sigma_3\mapsto\sigma_1, \quad \sigma_2\mapsto\sigma_2.  \eqno(\eqDefR)
$$
(we denote it by $R$ because, if we interpret $\B_n$ as $\pi_1(X_n)$ where
$X_n$ is the space of monic squarefree polynomials of degree $n$,
then $R$ is induced by the mapping which takes a degree $4$ polynomial
to its cubic resolvent).

For a group $G$, we denote its commutator subgroup, center, and abelianization
by $G'$, $Z(G)$, and $G^\ab$ respectively. We also denote
the inner automorphism $y\mapsto xyx^{-1}$ by $\tilde x$,
the commutator $xyx^{-1}y^{-1}$ by $[x,y]$, and the centralizer of 
an element $x$ (resp. of a subgroup $H$) in $G$
by $Z(x;G)$ (resp. by $Z(H;G)$).

Given two group homomorphisms $f:G_1\to G_2$ and
$\tau:G_1^\ab\to Z(\im f;G_2)$,
we define the {\it transvection} of $f$ by $\tau$ as the
homomorphism $f_{[\tau]}:G_1\to G_2$ given by $x\mapsto f(x)\tau(\bar x)$
where $\bar x$ is the image of $x$ in $G_1^\ab$. To simplify notation, we
will not distinguish between $\tau$ and its composition with the canonical
projection $G_1\to G_1^\ab$. So, we shall often speak of a transvection
by $\tau:G_1\to Z(\im f;G_2)$.

We say that two homomorphisms $f,g:G_1\to G_2$ are {\it equivalent}
if there exists $h\in\Aut(G_2)$ such that $f=hg$. If, moreover, $h\in\Inn(G_2)$,
we say that $f$ and $g$ are {\it conjugate}.

\proclaim{ Theorem \thBBfour }
Any homomorphism $\varphi:\B'_4\to\B_4$ either is equivalent to a transvection
of the inclusion map, or $\varphi=fR$ for a homomorphism $f:\B'_3\to\B_4$
(since $\B'_3$ is free [\refGL], it has plenty of homomorphisms to any group).
\endproclaim

We prove this theorem in \S\sectBfour.

\proclaim{ Corollary \corBBfour }
Any endomorphism of $\B'_4$ is either an automorphism or a composition
of $R$ with a homomorphism $\B'_3\to\B'_4$.
\endproclaim

As we already mentioned, $\B'_3$ is free, thus its homomorphisms are evident.
Now let us describe endomorphisms of $\B_n$.
We say that a homomorphism is {\it cyclic} if its image is a cyclic group
(probably, infinite cyclic).

\proclaim{ Theorem \thBn } {\rm(proven for $n\ge 6$ in [\refCastel]).} 
If $n\ge 5$, then any non-cyclic endomorphism of $\B_n$
is a transvection of an automorphism.
\endproclaim

For $n\ge 7$, this result is derived in [\refKM] from Theorem~\thBBn.
The same proof works without any change for any $n\ge 5$.

\proclaim{ Theorem \thBfour }
Any endomorphism of $\B_4$ is either a transvection of an automorphism, or
it is of the form $fR$ for some $f:\B_3\to\B_4$
(see Proposition~\propBthreeBfour\ for a general form of such $f$).
\endproclaim

This theorem also can be derived from Theorem~\thBBfour\
in the same way as it is done in [\refKM] for $n\ge 7$.

Let $\Delta=\Delta_n=\prod_{i=1}^{n-1}\prod_{j=1}^{n-i}\sigma_j$
(the Garside's half-twist),
$\delta=\delta_n=\sigma_{n-1}\dots\sigma_2\sigma_1$, and
$\gamma=\gamma_n=\sigma_1\delta_n$.
One has $\delta^n=\gamma^{n-1}=\Delta^2$, and it is known
that $Z(\B_n)$ is generated by $\Delta^2$, and each periodic
braid (i.e. a root of a central element) is conjugate to $\delta^k$ or
$\gamma^k$ for some $k\in\Z$.

It is well-known that $\B_3$ admits a presentation
$\langle\Delta,\delta\mid\Delta^2=\delta^3\rangle$.
By combining this fact with basic properties of canonical reduction systems,
it is easy to prove the following descriptions of homomorphisms from
$\B_3$ to $\B_n$ for $n=3$ or $4$.

\proclaim{ Proposition \propBthree }
Any non-cyclic endomorphism of $\B_3$ is equivalent to a transvection by $\tau$ of
a homomorphism of the form $\Delta\mapsto\Delta$, $\delta\mapsto X\delta X^{-1}$
for some $X\in\B_3$ and $\tau:\B_3^{\ab}\to Z(\B_3)=\langle \Delta^2\rangle$.
\endproclaim

\proclaim{ Proposition \propBthreeBfour }
For any non-cyclic homomorphism $\varphi:\B_3\to\B_4$, one of the following two
possibilities holds:
\roster
\item"(a)" $\varphi$ is equivalent to a transvection by $\tau$ of a homomorphism
of the form $\Delta_3\mapsto\Delta_4$, $\delta_3\mapsto X\gamma_4 X^{-1}$
for some $X\in\B_4$ and $\tau:\B_3^{\ab}\to Z(\B_4)=\langle \Delta_4^2\rangle$;

\item"(b)" $\varphi$ is equivalent to $(\iota\psi)_{[\tau]}$ 
where $\psi$ is a non-cyclic endomorphism of $\B_3$, $\iota:\B_3\to\B_4$ is the standard 
embedding, and $\tau$ is a homomorphism
$\B_3^{\ab}\to Z(\B_4) = \langle\Delta_4^2\rangle$.
\endroster
\endproclaim

\medskip\noindent
{\bf Remark \remTransv.}
Since $\B_n^\ab\cong Z(\B_n)\cong\Z$, the transvection
in Theorem~\thBn\ (and in the non-degenerate case in Theorem~\thBfour)
is uniquely determined by a single integer number.
In contrast, $(\B'_4)^\ab\cong\Z^2$, thus the transvection in Theorem~\thBBfour\
depends on two integers (here $\im\varphi=\B_4'$,
hence $Z(\im\varphi;\B_4)=Z(\B_4)\cong\Z$).
Notice also that two transvections are involved in the case (b) of
Proposition~\propBthreeBfour, thus the general form of $\varphi$ in this case is
$$
  \Delta_3\mapsto f\big(\iota(\Delta_3)^{6k+1}\Delta_4^{6l}\big),
  \qquad
  \delta_3\mapsto f\big(\iota(X\delta_3X^{-1}\Delta_3^{4k})\Delta_4^{4l}\big)
$$
with $k,l\in\Z$, $X\in\B_3$, $f\in\Aut(\B_4)$.

\medskip

\if{
When writing [\refAFST], I restricted myself to automorphisms of $\B'_n$
just because Vladimir Lin's question [\refLinTalk] that I answered to,
concerned the automorphisms only. By adding some more arguments, the same
proofs could be adapted to the classification of endomorphism.

However, a big amount of that proof is
no longer needed since 
Later on, Lin told me that the same question about endomorphisms of $\B'_n$ is also
interesting. Then I reread my paper and convinced myself that, by adding some
extra arguments, the same proof can be adapted for all endomorphisms as well.
However the paper was already written and submitted, and I did not do it
(of course, since the proof was not written, one cannot be sure that it was complete).
Then the problem was solved in [\refKM] for $n\ge 7$, and I decided
to write down a proof for smaller $n$. Since some intermediate results in [\refKM]
are valid for $n\ge 5$ (see Proposition~\propKM\ below),
using them I significantly reduced my proof in the case $n\ge 5$.
In contrast, the proof in the case $n=4$ remained the same, and even became
more complicated because, being inspired by [\refKM], I decided to extent it
from endomorphisms of $\B'_4$ to all homomorphisms $\B'_4\to\B_4$
(these are almost the same problems for $n\ge 5$ but not for $n=4$).

A specific feature of the problem of homomorphisms of $\B_n$ or $\B'_n$
is that for smaller $n$ it is somewhat harder than for larger $n$ because
the bigger $n$ is, the easier is to find commuting elements in the group.
So, for $n\ge 7$ a proof based almost uniquely on Nielsen--Thurston \
theory is possible, more
precisely, on the properties of canonical reduction systems of commuting
elements of mapping class groups (see [\refKM]).
For $n=6$, I needed to add some elementary linear algebra arguments (see Lemma~\lemJab).
I exposed them using the language of the representation theory of
symmetric groups but, of course, the same could be done just by solving certain
systems of linear equations, as it was done in [\refDG].
Another essential ingredient of the proof for $n=6$
are results of Lin [\refLin, Theorems~C and~D] on homomorphisms of pure
commutator subgroups.
For $n=5$, in addition to that, I did not manage to get rid of an argument coming
from the Garside theory (see Remark~\remGarside).
In the case $n=4$ there are so few pairs of commuting
elements, that a completely different approach is needed.
Fortunately, the lack of commuting pairs is
compensated by the semidirect product structure, which I extensively use in
my proof.
}\fi


\head\sectBn. The case $n\ge 5$
\endhead

In this section we prove Theorem~\thBBn\ which describes
homomorphisms $\B'_n\to\B_n$ for $n\ge 5$. The proof is very
similar to the proof  of the case $n\ge5$ of
the main theorem of [\refAFST] which describes $\Aut\B'_n$.
As we already mentioned, Theorem~\thBBn\ for $n\ge 7$ is proven
by Kordek and Margalit in [\refKM]. Some elements of their proof are
valid for $n\ge 5$ (see Proposition~\propKM\ below) which allowed us to
omit a big part of our original proof based on [\refAFST].

Let $\bS_n$ be the symmetric group. Let $e:\B_n\to\Z$ and
$\mu:\B_n\to\bS_n$ be the homomorphisms defined on the generators by
$e(\sigma_i)=1$ and $\mu(\sigma_i)=(i,i+1)$ for
$i=1,\dots,n-1$. So, $e(X)$ is the exponent sum (signed word length) of $X$.
Let $\P_n=\ker\mu$ be the pure braid group. Following [\refLin], we
denote $\P_n\cap\B'_n$ by $\J_n$, and $\mu|_{\B'_n}$ by $\mu'$, thus
$\J_n=\ker\mu'$.

For a pure braid $X$, we denote the linking number between the $i$-th and
the $j$-th strings of $X$ by $\lk_{ij}(X)$. It can be defined as
$\frac12e(X_{ij})$ where
$X_{ij}$ is the 2-braid obtained from $X$ by removal of all strings except
the $i$-th and the $j$-th ones.
For $1\le i<j\le n$, we set
$\sigma_{ij}=(\sigma_{j-1}\dots\sigma_{i+1})\sigma_i
(\sigma_{j-1}\dots\sigma_{i+1})^{-1}$ (here $\sigma_{i,i+1}=\sigma_i$). Then $\P_n$
is generated by $\{\sigma_{ij}^2\}_{1\le i<j\le n}$ (see [\refArtin])
and we denote the image of $\sigma_{ij}^2$ in $\P_n^\ab$ by $A_{ij}$.
We use the additive notation for $\P_n^\ab$ and $\J_n^\ab$.

\proclaim{ Lemma \lemPab } {\rm([\refAFST, Lemma~2.3]).}
$\P_n^\ab$ (for any $n$) is free abelian group with basis
$(A_{ij})_{1\le i<j\le n}$,
and the natural projection $\P_n\to\P_n^\ab$ is given by
$X\mapsto\sum_{i<j} \lk_{ij}(X)A_{ij}$.

If $n\ge 5$, then the homomorphism $\J_n^\ab\to\P_n^\ab$ induced by
the inclusion map defines an isomorphism of $\J_n^\ab$ with
$\{\sum x_{ij}A_{ij}\mid\sum x_{ij}=0\}$ {\rm(notice that this statement is wrong
for $n=3$ or $4$; see [\refAFST, Proposition~2.4])}.
\endproclaim

From now on, till the end of this section, we assume that
$n\ge 5$ and $\varphi:\B'_n\to\B_n$ is a non-cyclic homomorphism.
Since any group homomorphism $G_1\to G_2$ maps $G'_1$ to $G'_2$, we have
$\varphi(\B''_n)\subset\B'_n$. By [\refGL] (see also [\refAFST, Remark~2.2]),
we have $\B''_n=\B'_n$, thus
$$
     \varphi(\B'_n)\subset\B'_n.
$$
Then [\refLin, Theorem D] implies that
$$
     \varphi(\J_n)\subset\J_n.
$$
Thus we may consider the endomorphism $\varphi_*$ of $\J_n^\ab$ induced
by $\varphi|_{\J_n}$.
We shall not distinguish between $\J_n^\ab$ and its isomorphic
image in $\P_n^\ab$ (see Lemma~\lemPab).

Following [\refLin], we set
$$
        c_i = \sigma_1^{-1}\sigma_i
        \quad(i=3,\dots,n-1) \quad\text{and $c=c_3$.}
$$

\proclaim{ Lemma \lemJab }
Suppose that $\mu\varphi=\mu'$ and
$\varphi(c)=c$.
Then $\varphi_* = \id$.
\endproclaim

\demo{ Proof }
The exact sequence $1\to\J_n\to\B'_n\to\A_n\to1$
defines an action of $\A_n$ on $\J_n^\ab$ by conjugation.
 Let $V$ be a complex vector space with base
$e_1,\dots,e_n$ endowed with the natural action of $\bS_n$
induced by the action on the base. We identify $\P^\ab_n$ with
its image in the symmetric square $\Sym^2 V$ under the homomorphism
$A_{ij}\to e_i e_j$. Then, by Lemma~\lemPab, we may identify
$\J^\ab_n$ with $\big\{\sum x_{ij}e_ie_j\,\big|\,x_{ij}\in\Z, \sum x_{ij}=0\big\}$.
These identifications are compatible with the action of $\A_n$.
Thus $W:=\J_n^\ab\otimes\C$ is a $\C\A_n$-submodule of $\Sym^2 V$.

For an element $v$ of a $\C\bS_n$-module, let
$\langle v\rangle_{\C\bS_n}$ be the $\C\bS_n$-submodule generated by $v$.
It is shown in the proof of [\refAFST, Lemma 3.1],
that $W=W_2\oplus W_3$ where
$$
   W_2 = \langle(e_1-e_2)(e_3+\dots+e_n)\rangle_{\C\bS_n}, \qquad
   W_3 = \langle(e_1-e_2)(e_3-e_4)\rangle_{\C\bS_n},
$$
and that $W_2$ and $W_3$ 
are irreducible $\C\bS_n$-modules isomorphic to the
Specht modules corresponding to the partitions $(n-1,1)$ and $(n-2,2)$
respectively. Since the Young diagrams of these partitions are not symmetric,
$W_2$ and $W_3$ are also irreducible as $\C\A_n$-modules.

The condition $\mu\varphi=\mu'$ implies that $\varphi_*$ is
$\A_n$-equivariant.
Hence, by Schur's lemma, $\varphi_*=a\,\id_{W_2}\oplus\, b\,\id_{W_3}$.
We have the identity 
$$
  (n-2)(e_1-e_2)e_3=(e_1-e_2)(e_3+\dots+e_n)+\sum_{i\ge4}(e_1-e_2)(e_3-e_i)
$$
whence, denoting $e_5+\dots+e_n$ by $e$,
$$
\split
   (n-2)\varphi_*((e_1-e_3)e_2)&=
   (e_1-e_3)\big(a(e_2+e_4+e) + b((n-3)e_2-e_4-e)\big),
\\
   (n-2)\varphi_*((e_2-e_4)e_3)&=
   (e_2-e_4)\big(a(e_1+e_3+e) + b((n-3)e_3-e_1-e)\big).
\endsplit
$$
The condition $\varphi(c)=c$ implies
the $\varphi$-invariance of
$c^2\in\J_n$. Since the image of $c^{-2}$
in $\J_n^\ab$ is $A_{12}-A_{34}$, we obtain that $e_1e_2-e_3e_4$
is $\varphi_*$-invariant. Hence
$$
\split
 (n-2)(e_1e_2-e_3e_4)&=
 (n-2)\varphi_*(e_1e_2-e_3e_4)
\\&=
 (n-2)\varphi_*\big((e_1-e_3)e_2 + (e_2-e_4)e_3\big)
\\&=
 \big(2a+(n-4)b\big)(e_1e_2-e_3e_4) + (a-b)(e_1+e_2-e_3-e_4)e
\endsplit
$$
Since $\{e_ie_j\}_{i<j}$ is a base of $\Sym^2 V$, it follows that
$2a+(n-4)b=n-2$ and $a-b=0$ whence $a=b=1$.
\qed\enddemo

\proclaim{ Lemma \lemMu } Let $\varphi_1$ and $\varphi_2$ be
equivalent homomorphisms $\B'_n\to\B_n$. Then $\mu\varphi_1$ and
$\mu\varphi_2$ are conjugate.
\endproclaim

\demo{ Proof } This fact immediately follows from Dyer -- Grossman's [\refDG]
classification of automorphisms of $\B_n$ (see the beginning of the introduction)
because $\mu\Lambda=\mu$.
\qed\enddemo

\proclaim{ Proposition \propKM }
{\rm(Kordek and Margalit [\refKM, \S3, Proof of Thm.~1.1, Cases 1--3 and Step~1
of Case~4]).}
There exists $f\in\Aut(\B_n)$ such that
$f\varphi(c_i) = c_i$ for each
odd $i$ such that $1\le i<n$ (recall that we assume $n\ge 5$).
\endproclaim

This proposition implies, in particular, that
$\mu\varphi$ is non-trivial, hence by Lin's result
[\refLin, Theorem C]
$\mu\varphi$ is conjugate either to $\mu'$ or to $\nu\mu'$ (when $n=6$)
where $\nu$ is the restriction to $\A_6$ of the automorphism of $\bS_6$
given by $(12)\mapsto(12)(34)(56)$,
$(123456)\mapsto(123)(45)$ (it represents the only
nontrivial element of $\Out(\bS_6))$.

\proclaim{ Lemma \lemKMsix } If $n=6$, then $\mu\varphi$ is not conjugate to $\nu\mu'$.
\endproclaim

\demo{ Proof }
Let $H$ be the subgroup generated by $c_3$
and $c_5$. By Lemma~\lemMu\ and Proposition~\propKM\ we may
assume that $\varphi|_H=\id$.
Then we have
$$
  \mu'(H)=\mu\varphi(H)=\{\id,(12)(34),(12)(56),(34)(56)\}.
$$
In particular, no element of $\{1,\dots,6\}$ is fixed
by all elements of $\mu\varphi(H)$. A straightforward computation shows that
$$
  \nu\mu'(H)=\{\id,(12)(34),(13)(24),(14)(23)\},               \eqno(\eqKlein)
$$
thus $5$ and $6$ are fixed by all elements of $\nu\mu'(H)$. Hence these
subgroups are not conjugate in $\bS_6$.
\qed\enddemo

\proclaim{ Lemma \lemKM } There exists $f\in\Aut(\B_n)$ such that
$f\varphi(c)=c$ and $\mu f\varphi=\mu'$.
\endproclaim

\demo{ Proof } By Proposition~\propKM\ we may assume that
$$
     \varphi(c)=c.                                      \eqno(\eqKM)
$$
Then $\mu\varphi$ is non-trivial,
hence, by [\refLin, Thm.~C] combined with Lemma~\lemKMsix,
it is conjugate to $\mu'$, i.e. there exists
$\pi\in\bS_n$ such that $\tilde\pi\mu\varphi=\mu'$, i.e.
$\pi\mu(\varphi(x)) = \mu(x)\pi$ for each $x\in\B'_n$. For
$x=c$ this implies by (\eqKM) that $\pi$ commutes with $(12)(34)$,
hence $\pi=\pi_1\pi_2$ where $\pi_1\in V_4$ (the group in the right hand size of
(\eqKlein)) and $\pi_2(i)=i$ for $i\in\{1,2,3,4\}$.
Let $\tilde V_4=\{1,c,\Delta_4,c\Delta_4\}$. This is not a subgroup but
we have $\mu(\tilde V_4)=V_4$.
We can choose $y_1\in\tilde V_4$ and $y_2\in\langle\sigma_5,\dots,\sigma_{n-1}\rangle$
so that $\mu(y_j)=\pi_j$, $j=1,2$. Let $y=y_1y_2$.
Then we have $\tilde y(c)=c^{\pm 1}$ and $\mu\tilde y\varphi=\tilde\pi\mu\varphi=\mu'$.
Thus, for $f=\Lambda^k\tilde y$, $k\in\{0,1\}$, we have
$f\varphi(c)=c$ and $\mu f\varphi=\mu'$.
\qed\enddemo

Due to Lemma~\lemKM, from now on we assume that $\mu\varphi=\mu'$ and $\varphi(c)=c$.
Then, by Lemma~\lemJab, we have $\varphi_*=\id$, hence (see Lemma~\lemPab)
$$
      \lk_{ij}(x) = \lk_{ij}(\varphi(x))
      \qquad\text{for any $x\in\J_n$ and $1\le i<j\le n$.}           \eqno(\eqLK)
$$

Starting at this point, the proof of [\refAFST, Thm.~1.1]
given in [\refAFST, \S5], can be repeated almost word-by-word in our setting.
The only exception is the proof of [\refAFST, Lemma~5.8]
(which is Lemma~\lemEight\ below) where
the invariance of the isomorphism type of centralizers of certain elements
is used as well as Dyer--Grossman result [\refDG].
However, as pointed out in [\refAFST, Remark~5.15]
(there is a misprint there: $n\ge 6$ should be replaced by $n\ge 5$),
there is another, even simpler, proof of Lemma~\lemEight\
based on Lemma~\lemXYX\ (see below).
This proof was not included in [\refAFST] by the
following reason. At that time we new only Garside-theoretic proof of
Lemma~\lemXYX\ while the rest of the proof of the main theorem
for $n\ge 6$ used only Nielsen-Thurston theory and results of [\refLin].
So we wanted to make
the proofs (at least for $n\ge6$) better accessible for readers who are not familiar
with the Garside theory. Now we learned from [\refKM] that when we wrote that paper,
Lemma~\lemXYX\ had been already known for a rather long time [\refBM, Lemma~4.9]
and the proof in [\refBM] is based on Nielsen-Thurston theory. 

In the rest of this section, for the reader's convenience we re-expose 
Section 5.1 of [\refAFST] (Sections~5.2--5.3 can be left without any change).
In this re-exposition we give another proof of [\refAFST, Lemma~5.8]
and omit the lemmas which are no longer needed due to Proposition~\propKM.

We shall consider $\B_n$ as a mapping class group of $n$-punctured disk $\Bbb D$.
We assume that $\Bbb D$ is a round disk in $\C$ and the set of the punctures
is $\{1,2,\dots,n\}$.
Given an embedded segment $I$ in $\Bbb D$ with endpoints at two punctures, we
denote with $\sigma_I$ the positive half-twist along the boundary of
a small neighborhood of $I$. The set of all such braids is the conjugacy
class of $\sigma_1$ in $\B_n$.
The arguments in the rest of this section
are based on Nielsen-Thurston theory. The main tool are the
canonical reduction systems. One can use [\refBLM], [\refFM] or [\refI]
as a general introduction to the subject. In [\refAFST] we gave all
precise definitions and statements needed there (using the language and notation
inspired mostly by [\refGMW]).

\proclaim{ Lemma~\lemXYX } {\rm([\refBM, Lemma~4.9], [\refAFST, Lemma~A.2]).}
Let $x,y\in\B_n$ be such that $xyx=yxy$ and each of $x$ and $y$ is conjugate to $\sigma_1$.
Then there exists $u\in\B_n$ such that $\tilde u(x)=\sigma_1$ and $\tilde u(y)=\sigma_2$.
\endproclaim



Let $\sh_2:\B_{n-2}\to\B_n$ be the homomorphism
$\sh_2(\sigma_i)=\sigma_{i+2}$. We set
$$
     \tau=\sigma_1^{(n-2)(n-3)}\sh_2(\Delta_{n-2}^{-2}).
$$
We have $\tau\in\J_n$
(in the notation of [\refAFST],
$\tau=\psi_{2,n-2}(1;\sigma_1^{(n-2)(n-3)},\Delta^{-2})$).
Recall that we assume $\varphi(c)=c$, $\mu\varphi=\mu'$, and hence
(\eqLK) holds.

\proclaim{ Lemma~\lemOne }
Let $I$ and $J$ be two disjoint embedded segments with endpoints at punctures.
Then $\varphi(\sigma_I^{-1}\sigma_J)=\sigma_{I_1}^{-1}\sigma_{J_1}$ where
$I_1$ and $J_1$ are disjoint embedded segments such that
$\partial I_1=\partial I$ and $\partial J_1=\partial J$.
\endproclaim

\demo{ Proof } The braid $\sigma_I^{-1}\sigma_J$ is conjugate to $c$,
hence so is its image (because $\varphi(c)=c)$. Therefore
$\varphi(\sigma_I^{-1}\sigma_J)=\sigma_{I_1}^{-1}\sigma_{J_1}$ for
some disjoint $I_1$ and $J_1$. The matching of the boundaries follows from
(\eqLK) applied to $\sigma_I^{-2}\sigma_J^2$.
\qed\enddemo

\proclaim{ Lemma~\lemThree } {\rm(cf.~[\refAFST, Lemmas~5.1 and 5.3]).}
Let $C_1$ be a component of the canonical reduction system of $\varphi(\tau)$.
Then $C_1$ cannot separate the punctures $1$ and $2$, and it
cannot separate the punctures $i$ and $j$ for  $3\le i<j<n$.
\endproclaim

\demo{ Proof }
Let $u=\sigma_1^{-1}\sigma_{ij}$, $3\le i<j\le n$.
By Lemma~\lemOne, $\varphi(u)=\sigma_I^{-1}\sigma_J$ with
$\partial I=\{1,2\}$ and $\partial J=\{i,j\}$. Since
$\varphi(u)$ commutes with $\varphi(\tau)$, the result follows.
\qed\enddemo

\proclaim{ Lemma~\lemSeven } {\rm(cf.~[\refAFST, Lemma~5.7]).}
$\varphi(\tau)$ is conjugate in $\P_n$ to $\tau$.
\endproclaim

\demo{ Proof }
$\varphi(\tau)$ cannot be pseudo-Anosov because
it commutes with $\varphi(c)$ which is $c$ by our assumption,
hence it is reducible.

If $\varphi(\tau)$ were periodic, then it would be
a power of $\Delta^2$ because it is a pure braid.
This contradicts (\eqLK), hence $\varphi(\tau)$ is reducible non-periodic.

Let $C$ be the canonical reduction system for $\varphi(\tau)$.
By Lemma \lemThree, one of the following three cases occurs.

\smallskip
Case 1. $C$ is connected, the punctures $1$ and $2$ are inside $C$, all the other punctures are outside $C$.
Then the restriction of $\varphi(\tau)$ (viewed as a diffeomorphism of $\Bbb D$)
to the exterior of $C$
cannot be pseudo-Anosov because $\varphi(\tau)$ commutes with
$\varphi(c)=c$, hence it preserves a circle which separates $3$ and $4$ from
$5,\dots,n$. Hence $\varphi(\tau)$ is periodic which contradicts (\eqLK).
Thus this case is impossible.

\smallskip
Case 2. $C$ is connected, the punctures $1$ and $2$ are outside $C$, all the other punctures are inside $C$.
This case is also impossible and the proof is almost the same as in Case 1.
To show that $\varphi(\tau)$ cannot be pseudo-Anosov,
we note that it preserves a curve which encircles only $1$ and $2$.

\smallskip
Case 3. $C$ has two components: $C_1$ and $C_2$ which encircle $\{1,2\}$ and
$\{3,\dots,n\}$ respectively.
Let $\alpha$ be the interior braid of $C_2$ (that is $\varphi(\tau)$
with the strings $1$ and $2$ removed).
It cannot be pseudo-Anosov by the same reasons as in Case 1:
because $\varphi(\tau)$ preserves a circle separating $3$ and $4$ from $5,\dots,n$.
Hence $\alpha$ is periodic. Using (\eqLK),
we conclude that $\varphi(\tau)$ is a conjugate of $\tau$.
Since the elements of $Z(\tau;\B_n)$ realize any permutation of $\{1,2\}$ and
of $\{3,\dots,n\}$, the conjugating element can be chosen in $\P_n$.
\qed\enddemo

\proclaim{ Lemma~\lemEight } {\rm(cf.~[\refAFST, Lemma~5.8]).}
There exists $u\in\P_n$ such that
$\varphi(c_i)=\tilde u(c_i)$ for
each $i=3,\dots,n-1$.
\endproclaim

\demo{ Proof }
Due to Lemma \lemSeven, without loss of generality we may assume that
$\varphi(\tau)=\tau$ and $\tau(C)=C$ where $C$ is the canonical reduction system
for $\tau$ consisting of two round circles
$C_1$ and $C_2$ which encircle $\{1,2\}$ and $\{3,\dots,n\}$ respectively.
Since the conjugating element in Lemma~\lemSeven\ is chosen in $\P_n$, we
may assume that (\eqLK) still holds.

By Lemma~\lemOne, for each $i=3,\dots,n-1$, we have
$\varphi(c_i)=\sigma_{I_i}^{-1}\sigma_{J_i}^{-1}$ with
$\partial I_i=\{1,2\}$ and $\partial J_i=\{i,i+1\}$.
Since $\tau$ commutes with each $c_i$, the segments $I_i$ and $J_i$
can be chosen disjoint from the circles $C_1$ and $C_2$. Hence
$\sigma_{I_i}=\sigma_1$ for each $i$, and all the segments $J_i$ are
inside $C_2$.

Therefore the braids $\sigma_{J_3},\dots,\sigma_{J_{n-1}}$ satisfy the
same braid relations as $\sigma_3,\dots,\sigma_{n-1}$. Hence, by
Lemma~\lemXYX\ combined with [\refAFST, Lemma~5.13], 
$J_3\cup\dots\cup J_{n-1}$ is an embedded segment.
Hence it can be transformed to the straight line segment $[3,n]$
by a diffeomorphism identical on the exterior of $C_2$. Hence for the
braid $u$ represented by this diffeomorphism we have $\tilde u(c_i)=c_i$, $i\ge 3$.
The condition $\partial J_i=\{i,i+1\}$ implies that $u\in\P_n$.
\qed\enddemo

The rest of the proof of Theorem~\thBBn\ repeats word-by-word [\refAFST, \S\S5.2--5.3].

\medskip\noindent
{\bf Remark~\remGarside.} Besides Nielsen-Thurston theory, in the case $n=5$, the
arguments in [\refAFST, \S5.3] use an auxiliary result
[\refAFST,~Lemma~A.1] for which the only proof we know is based on
a slight modification of the main theorem of [\refJAlg] which is proven there
using the Garside theory.


\medskip



\head\sectBfour. The case $n=4$
\endhead

We shall use the same notation as in [\refAFST, \S6].
The groups $\B'_3$ and $\B'_4$ were computed in [\refGL], namely
$\B'_3$ is freely generated by
$u=\sigma_2\sigma_1^{-1}$ and $t=\sigma_1^{-1}\sigma_2$, and
$\B'_4=\bold K_4\rtimes\B'_3$ where
$\bold K_4=\ker R$ (see (\eqDefR)).
The group $\bold K_4$ is freely generated by
$c=\sigma_3\sigma_1^{-1}$ and $w=\sigma_2c\,\sigma_2^{-1}$.
The action of $\B'_3$ on $\bold K_4$ by conjugation is given by
$$
  ucu^{-1}=w,          \qquad
  uwu^{-1}=w^2c^{-1}w, \qquad
  tct^{-1}=cw,         \qquad
  twt^{-1}=cw^2.                                      \eqno(\eqGL)
$$
The action of $\sigma_1$ and $\sigma_2$ on $\K_4$ is given by
$$
  \sigma_1 c \sigma_1^{-1} = c,         \quad
  \sigma_1 w \sigma_1^{-1} = c^{-1}w,   \quad
  \sigma_2 c \sigma_2^{-1} = w,         \quad
  \sigma_2 w \sigma_2^{-1} = wc^{-1}w.                \eqno(\eqGLsig)
$$
So, we also have $\B_4=\K_4\rtimes\B_3$.

Besides the elements $c,w,u,t$ of $\B'_4$, we consider also
$$
   d
    =\Delta\sigma_1^{-3}\sigma_3^{-3}
   \qquad\text{and}\qquad
   g=R(d)
    =\Delta_3^2\sigma_1^{-6}.\;\;
$$
\if01{
The {\it mixed braid groups} $\B_{\vec m}$, $\vec m=\{m_1,\dots,m_k\}$
(which are subgroups of $\B_{m_1+\dots+m_k}$) and the {\it cabling maps}
$$
   \psi_{\vec m}:\B_k\times\B_{m_1}\times\dots\times\B_{m_k}\to
   \B_{\vec m}
$$
(which are injective) are defined in [\refAFST, \S2.3].
Here we shall use only $\B_{2,1,1}$ and $\B_{2,2}$
}\fi
One has (see Figure~\figDTU)
$$
     d=[c^{-1}t,u^{-1}],\qquad g=[t,u^{-1}].           \eqno(\eqDTU)
$$
We denote the subgroup generated by $c$ and $d$ by $H$ and
the subgroup generated by $c$ and $g$ by $G$.

\midinsert
\centerline{\epsfxsize=40mm\epsfbox{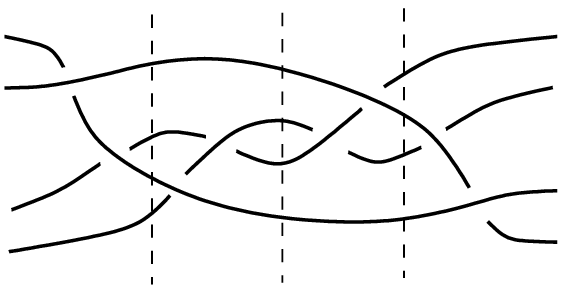}}
\vskip-3pt
\centerline{   \hskip3mm $c^{-1}t$
               \hskip4mm $u^{-1}$
               \hskip2mm $t^{-1}c$
               \hskip4mm $u$ \hskip7mm}
\botcaption{Figure \figDTU}  the identity $d=[c^{-1}t,u^{-1}]$.
\endcaption
\endinsert

Let $\varphi:\B'_4\to\B_4$ be a homomorphism such that
$\K_4\not\subset\ker\varphi$.

\proclaim{ Lemma \lemInj }
The restriction of $\varphi$ to $H$ 
is injective, $\varphi(H)\subset\B'_4$, and $\varphi(G)\subset\B'_4$.
\endproclaim

\demo{ Proof }
We have $H=\langle c\rangle\rtimes\langle d\rangle$ and $d$ acts on $c$
by $dcd^{-1}=c^{-1}$. Hence
any non-trivial normal subgroup of $H$
contains a power of $c$. Thus, if
$\varphi|_H$ 
were not injective, $\ker\varphi$ would contain a power of $c$ and hence
$c$ itself because the target group $\B_4$
does not have elements of finite order.
Then we also have $w\in\ker\varphi$ because $w=ucu^{-1}$.
This contradicts the assumption $\K_4=\langle c,w\rangle\not\subset\ker\varphi$,
thus $\varphi|_H$ is injective.

We have $dcd^{-1}=c^{-1}$, hence the image of $\varphi(c)$ under the abelianization
$e:\B_4\to\Z$ is zero, i.e., $\varphi(c)\in\B'_4$. By (\eqDTU) we also have
$\varphi(d)\in\B'_4$ and $\varphi(g)\in\B'_4$, thus
$\varphi(H)\subset\B'_4$ and $\varphi(G)\subset\B'_4$.
\qed\enddemo

\proclaim{ Lemma \lemCommCG } $\varphi(c)$ and $\varphi(g)$ do not commute.
\endproclaim

\demo{ Proof } Suppose that $\varphi(c)$ and $\varphi(g)$ commute.
Then $\varphi(c)=\varphi(gcg^{-1})$. Hence (see Figure~\figCGCC)
$\varphi(c)=\varphi(w^{-1}c^{-1}w)$, i.e., $\varphi$ factors through
the quotient of $\B'_4$ by the relation $wc = c^{-1}w$. Let us denote
this quotient group by $\hat\B'_4$.

\midinsert
\centerline{\epsfxsize110mm\epsfbox{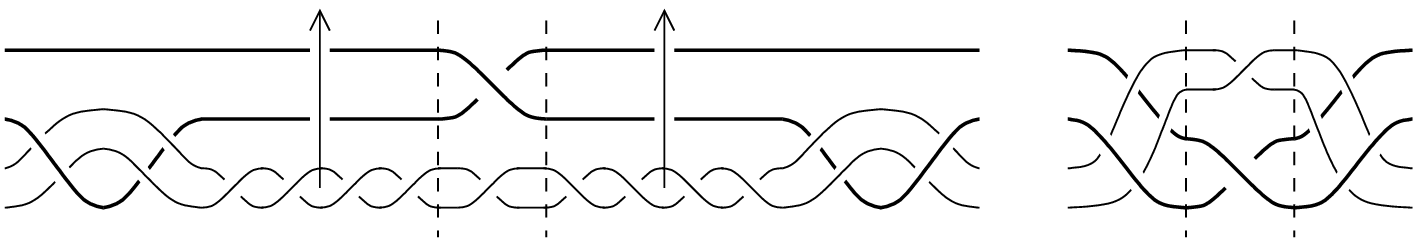}}
\vskip-3pt
\centerline{\hskip22mm $g$\hskip10mm $c$ \hskip10mm $g^{-1}$
            \hskip27mm $w^{-1}$ \hskip 2mm $c^{-1}$ \hskip3mm $w$}
\botcaption{Figure \figCGCC} The identity $gcg^{-1}=w^{-1}c^{-1}w$.
\endcaption
\endinsert

The relation $wc = c^{-1}w$ allows us to put any word
$\prod_j c^{k_j}w^{l_j}$ with $l_j=\pm1$
into the normal form $c^{k_1-k_2+k_3-\dots}w^{l_1+l_2+l_3+\dots}$
in $\hat\B'_4$.
Due to (\eqGL), the conjugation by $t$ of the word $w^{-1}cwc$
(which is equal to $1$ in $\hat\B'_4$) yields
$$
    1 = t(w^{-1}cwc)t^{-1} = (w^{-2}c^{-1})(cw)(cw^2)(cw) = w^{-1}cw^2cw=c^{-2}w^2
$$
(here in the last step we put the word into the above normal form).
Conjugating once more by $t$
and putting the result into the normal form, we get
$$
    1 = t(c^{-2}w^2)t^{-1} = (w^{-1}c^{-1})(w^{-1}c^{-1})(cw^2)(cw^2)
      =w^{-1}c^{-1}wcw^2 = c^2w^2.
$$
Thus $c^{-2}w^2 = c^2w^2 = 1$, i.e., $c^4=1$ in $\hat\B'_4$, hence
$\varphi(c^4)=1$ which contradicts Lemma~\lemInj.
\qed\enddemo

As in [\refAFST], we denote the stabilizer of $1$
under the natural action of $\B_3$ on $\{1,2,3\}$ by $\B_{1,2}$.
It is well-known (and easy to prove by Reidemeister-Schreier method) that
$\B_{1,2}$ is isomorphic to the Artin group of type $B_2$, that is
$\langle x,y\mid xyxy=yxyx\rangle$. The Artin generators $x$ and $y$
of the latter group correspond to $\sigma_1^2$ and $\sigma_2$.

\proclaim{ Lemma \lemG } {\rm(cf.~[\refAFST, Lemma 6.2])}
We have $G=Z(d^2c^6;\B'_4)$ and this group is generated by $g$ and $c$
subject to the defining relation $gcgc=cgcg$.
\endproclaim

\demo{ Proof }
The centralizer of $d^2c^6$ in $\B_4$ is the stabilizer of its canonical
reduction system which is shown in Figure~\figCRS, and
(see [\refGMW, Thm.~5.10]) it is the image of the injective homomorphism
$\B_{1,2}\times\Z\to\B_4$,
$(X,n)\mapsto Y\sigma_1^n$, where the \hbox{4-braid} $Y$ is obtained from the
\hbox{3-braid} $X$ by doubling the first strand.
It follows that $Z(d^2c^6;\B'_4)$ is the isomorphic image of $\B_{1,2}$ under
the homomorphism $\psi:\B_{1,2}\to\B'_4$ defined on the generators by
$\psi(\sigma_1^2)=g$, $\psi(\sigma_2)=c$ (see Figure~\figCG), thus $Z(d^2c^6;\B'_4)=G$.
As we have pointed out above, $\B_{1,2}$ is the Artin group of type $B_2$, hence
so is $G$ and $gcgc=cgcg$ is its defining relation.
\qed\enddemo

\midinsert
\centerline{\epsfxsize56mm\epsfbox{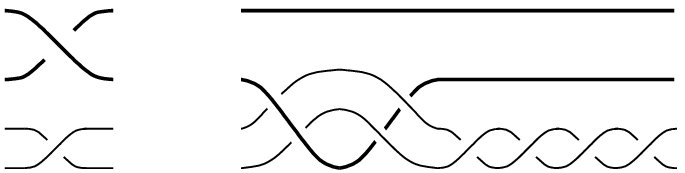}}
\centerline{$\psi(\sigma_2)=c$ \hskip16mm $\psi(\sigma_1^2)=g$ \hskip10mm}
\botcaption{Figure \figCG}
 The images of the generators under
                  $\psi:\B_{1,2}\to\B'_4$.
\endcaption
\endinsert

\proclaim{ Lemma \lemCD }
$\varphi(d^2c^6)$ is conjugate in $\B_4$ to $d^{2k}$, $d^{2k}c^{6k}$,
or $h^k$ for
some integer $k\ne0$,
where $h=\Delta_3^{-2}\sigma_3\sigma_2\sigma_1^2\sigma_2\sigma_3$.
\endproclaim

\demo{ Proof } Let $x=d^2c^6$. By Lemma~\lemG, $G=Z(x;\B'_4)$,
hence $\varphi(G)\subset Z(\varphi(x);\B_4)$.
By Lemma~\lemInj\ we also have $\varphi(G)\subset\B'_4$,
hence $\varphi(G)\subset Z(\varphi(x);\B'_4)$.
Then it follows from Lemma~\lemCommCG\ that
$Z(\varphi(x);\B'_4)$ is non-commutative.
The isomorphism classes of the centralizers (in $\B'_4$) of all elements
of $\B'_4$ are computed in [\refAFST, Table~6.1]. We see in this table that
$Z(\varphi(x);\B'_4)$ is non-commutative only in the required cases
(see the corresponding canonical reduction systems in Figure~\figCRS) 
unless $\varphi(x)=1$.
However the latter case is impossible by Lemma~\lemInj.
\qed\enddemo

\midinsert
\centerline{
   \epsfxsize=25mm\epsfbox{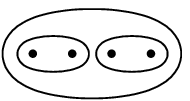}
   \hskip 1cm
   \epsfxsize=25mm\epsfbox{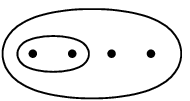}
   \hskip 1cm
   \epsfxsize=25mm\epsfbox{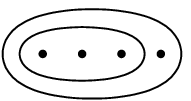}
}
\centerline{$d^{m}$ and $c^m$ \hskip22mm $(d^2c^6)^m$ \hskip26mm $h^m$ \;\;\;\;}
\botcaption{ Figure~\figCRS } Canonical reduc.~systems for
              $d^{m}$, $c^m$, $(d^2c^6)^m$, $h^m$, $m\ne0$.
\endcaption
\endinsert

\proclaim{ Lemma \lemCDii }
There exists an automorphism of $\B_4$ which takes $\varphi(c)$
and $\varphi(d)$ to $c^k$ and $d^k$ respectively for an odd positive integer $k$.
\endproclaim

\demo{ Proof }
Let $x=d^2c^6$ and $y=d^2c^{-6}$.
Since $y=dxd^{-1}$, the images of $x$ and $y$ are conjugate and
both of them belong to one of the conjugacy classes indicated in Lemma~\lemCD.
The canonical reduction systems for $d^{2k}$, $d^{2k}c^{6k}$, and
$h^k$ for $k\ne 0$ are
shown in Figure~\figCRS. Since $x$ and $y$ commute, the canonical
reduction systems of their images can be chosen disjoint from each other.
Hence, up to composing $\varphi$ with an inner automorphism of $\B_4$,
$\big(\varphi(x),\varphi(y)\big)$ is either $(h^{k_1},h^{k_2})$ or
$\big(d^{2k_1}c^{l_1},d^{2k_2}c^{l_2}\big)$ where $l_j\in\{0,\pm 6k_j\}$, $j=1,2$.
Since $x$ and $y$ are conjugate, by comparing the linking numbers
between different pairs of strings, we deduce that
$k_1=k_2$ and (in the second case) $l_1=\pm l_2$.
Moreover, $\varphi(x)\ne\varphi(y)$ by Lemma~\lemInj.
Hence, up to exchange of $x$ and $y$ (which is realizable by composing
$\varphi$ with $\tilde d$), we have
$\varphi(x)=d^{2k}c^{6k}$ and $\varphi(y)=d^{2k}c^{-6k}$ whence,
using that $xy^{-1}=c^{12}$, we obtain
$\varphi(c^{12})=\varphi(xy^{-1}) = c^{12k}$,
Since the canonical reduction systems of any braid and its non-zero power
coincide (see, e.g., [\refGM, Lemmas~2.1--2.3]), we obtain
$\varphi(c)=c^k$ and $\varphi(d)=d^k$. By composing $\varphi$ with $\Lambda$
if necessary, we can arrive to $k>0$.
The relation $d^kc^kd^{-k}=c^{-k}$ combined with Lemma~\lemInj\
implies that $k$ is odd.
\qed\enddemo

\proclaim{ Lemma \lemK }
 $\varphi(\K_4)\subset \K_4$.
\endproclaim

\demo{ Proof }
Lemma~\lemCDii\ implies that $c^k$ is mapped to $\varphi(c)$ by an automorphism of
$\B_4$. Since $\K_4$ is a characteristic subgroup of $\B'_4$
(see [\refAFST, Lemma 6.5]) and $\B'_4$ is a characteristic subgroup of $\B_4$,
we deduce that $\varphi(c)\in\K_4$. The same arguments can be applied to any other
homomorphism of $\B'_4$ to $\B_4$ whose kernel does not contain $\K_4$,
in particular, they can be applied to $\varphi\tilde u$
whence $\varphi\tilde u(c)\in\K_4$. Since $\varphi(w)=\varphi\tilde u(c)$,
we conclude that
$\varphi(\K_4)=\langle\varphi(c),\varphi(w)\rangle\subset\K_4$.
\qed\enddemo




Let
$$
                   F=G\cap\K_4.
$$

\proclaim{ Lemma \lemF }
(a). The group $F$ is freely generated by $c$ and $c_1=w^{-1}c^{-1}w$.

\smallskip
(b). Let $a_1,\dots,a_{m-1}$ and $b_1,\dots,b_m$ be non-zero integers, and
let $a_0$ and $a_m$ be any integers.
Then $c^{a_0}w^{b_1}c^{a_1}\dots w^{b_m}c^{a_m}$ is in $F$ if and only if
$m$ is even and $b_j=(-1)^j$ for each $j=1,\dots,m$.
\endproclaim

\demo{ Proof }
The relation on $g$ and $c$ in Lemma~\lemG\ is equivalent to
$$
    g^{-1}cgc = cgcg^{-1}.                         \eqno(\eqCGC)
$$
Recall that $G=\langle c,g\rangle$. We have $R(c)=1$ and,
by (\eqDTU), $g=R(d)\in\B'_3$ whence $R(g)=g$.
Hence $R(G)$ is generated by $g$.
By definition, $F=\ker(R|_G)$,
hence $F$ is the normal closure of $c$ in $G$, i.e., $F$ is
generated by the elements $\tilde g^k(c)$, $k\in\Z$.
We have $\tilde g(c) = c_1$ (see Figure~\figCGCC) and
$$
  \tilde g(c_1)=\tilde g^2(c)=
             g\,c^{-1}(cgcg^{-1})g^{-1} \overset{\text{by (\eqCGC)}}\to= 
             g\,c^{-1}(g^{-1}cgc)g^{-1} = 
              c_1^{-1}c\,c_1
$$
whence by induction we obtain
$\tilde g^k(c)\in\langle c,c_1\rangle$ for all positive $k$.
Similarly, 
$$
 \tilde g^{-1}(c)=(g^{-1}cgc)c^{-1}
                 \overset{\text{by (\eqCGC)}}\to=(cgcg^{-1})c^{-1}
                 =c(gcg^{-1})c^{-1}
                 =c\,c_1c^{-1}
$$
and $\tilde g^{-1}(c_1)=c$
whence
$\tilde g^k(c)\in\langle c,c_1\rangle$ for all negative $k$.
Thus $F=\langle c,c_1\rangle$.

To check that $c$ and $c_1$ is a free base of $F$ (which completes the proof of (a)),
it is enough to observe that if, in a reduced word in $x$, $y$, we replace
each $x^k$ with $c^k$ and each $y^k$ with $w^{-1}c^{-k} w$,
then we obtain a reduced word in $c$ and $w$. The statement (b) also easily
follows from this observation.
\qed\enddemo

\proclaim{ Lemma \lemXF } If $x\in F$ and $x=[w^{-1},A]$ with $A\in\K_4$, then
$x=[w^{-1},c^k]$, $k\in\Z$.
\endproclaim

\demo{ Proof } Let $A=w^{b_1}c^{a_1}\dots w^{b_m}c^{a_m}w^{b_{m+1}}$, $m\ge 0$,
where $a_1,\dots,a_m$ and $b_2,\dots,b_m$ are non-zero while
$b_1$ and $b_{m+1}$ may or may not be zero. If $m=0$, then $[w^{-1},A]=1=[w^{-1},c^0]$
and we are done.
If $m=1$, then $[w^{-1},A]=w^{b_1-1}c^{a_1}w\,c^{-a_1}w^{-b_1}$ where, by Lemma~\lemF(b),
we must have $b_1=0$, hence $[w^{-1},A]=[w^{-1},c^{a_1}]$ as required.
Suppose that $m\ge 2$. Then
$$
  [w^{-1},A]=w^{b_1-1}c^{a_1}\dots w^{b_m}c^{a_m} w\,
     c^{-a_m} w^{-b_m}\dots c^{-a_1}w^{-b_1}
$$
and this is a reduced word in $c$, $w$. Hence, by Lemma~\lemF(b), the sequence
of the exponents of $w$ in this word
(starting form $b_1-1$ when $b_1\ne1$ or from $b_2$ when $b_1=1$)
should be $(-1,1,-1,1,\dots,-1,1)$.
Such a sequence cannot contain $(\dots,b_m,1,-b_m,\dots)$. A contradiction.
\qed\enddemo

\proclaim{ Lemma \lemWF }
If $\varphi(d^2)=d^2$ and $\varphi(c)=c$, then $w^{-1}\varphi(w)\in F$.
\endproclaim

\demo{ Proof } For any $k\in\Z$ we have
$$
   \sigma_3^k w
     =\sigma_3^k(\sigma_2\sigma_3)(\sigma_1^{-1}\sigma_2^{-1})
     =(\sigma_2\sigma_3)\sigma_2^k(\sigma_1^{-1}\sigma_2^{-1})
     =(\sigma_2\sigma_3)(\sigma_1^{-1}\sigma_2^{-1})\sigma_1^k
     =w\sigma_1^k,
$$
hence $\sigma_3^k w \sigma_1^{-k}=w=\sigma_3^{-k}w\sigma_1^k$ and we obtain
$$
   d^2 w d^{-2} = \Delta^2\sigma_1^{-6}(\sigma_3^{-6}w\sigma_1^6)\sigma_3^6\Delta^{-2}
   = \sigma_1^{-6}(\sigma_3^{6}w\sigma_1^{-6})\sigma_3^6
   = c^6 w c^6.                                                     \eqno(\eqCWC)
$$
Set $x=w^{-1}\varphi(w)$, i.e., $\varphi(w)=wx$.
The relation (\eqCWC) combined with our hypothesis on $c$ and $d^2$
implies
$$
  c^6 wx c^6 = \varphi(c^6wc^6)=
   \varphi(\tilde d^2(w))=\tilde d^2(wx)
   =\tilde d^2(w)\tilde d^2(x)=c^6wc^6 d^2xd^{-2}
$$
whence $x(c^6d^2)=(c^6d^2)x$, i.e., $x\in Z(d^2c^6)$.
On the other hand, $\varphi(w)\in\K_4$ by Lemma~\lemK, hence
$x=w^{-1}\varphi(w)\in\K_4$. By Lemma~\lemG\ we have $Z(d^2;\B'_4)=G$,
thus $x\in Z(d^2c^6)\cap\K_4=G\cap\K_4=F$.
\qed\enddemo

\proclaim{ Lemma \lemUT } There exists $f\in\Aut(\B_4)$ and
a homomorphism $\tau:\B'_4\to Z(\B_4)$
such that $f\varphi(c)=c$, $f\varphi(d^2)=d^2$, and
$Rf\varphi=R\id_{[\tau]}$.
\endproclaim

\demo{ Proof }
By Lemma~\lemCDii\ we may assume that $\varphi(c)=c^k$ and $\varphi(d)=d^k$
for an odd positive $k$.
For $x\in\bold K_4$, we denote its image in $\bold K_4^\ab$
by $\bar x$ 
and we use the additive notation for $\bold K_4^\ab$.
Consider the homomorphism $\pi:\B_4\to\Aut(\bold K_4^\ab)=\GL(2,\Z)$,
where $\pi(x)$ is defined as the automorphism of $\K_4^\ab$ induced by
$\tilde x$;
here we identify  $\Aut(\bold K_4^\ab)$ with $\GL(2,\Z)$
by choosing $\bar c$ and $\bar w$ as a base of $\bold K_4^\ab$.
By Lemma~\lemK, $\varphi(w)\subset\K_4$, hence we may write
$\overline{\varphi(w)}=p\bar c+q\bar w$ with $p,q\in\Z$.
Then, for any $x\in\B_4$, we have
$$
     \pi\varphi(x).P = P.\pi(x) \qquad\text{where}\quad
     P=\left(\matrix k&p\\0&q\endmatrix\right).               \eqno(\eqPXP)
$$
($P$ is the matrix of the endomorphism of $\K_4^\ab$ induced by $\varphi|_{\K_4}$).
By (\eqCWC) we have
$$
  \pi(d^2)= \left(\matrix 1 & 12\\
                          0 &  1\endmatrix\right) \qquad\text{hence}\qquad
    \pi(d^{2k}).P - P.\pi(d^2) =
    \left(\matrix 0 & 12k(q-1)\\0 & 0\endmatrix\right).       \eqno(\eqPGP)
$$
Since $\varphi(d^2)=d^{2k}$,
we obtain from (\eqPXP) combined with (\eqPGP) that $q=1$,
i.e., $\overline{\varphi(w)}=p\bar c +\bar w$.
By (\eqGL) we have $\varphi(u)c^k\varphi(u)^{-1}=\varphi(ucu^{-1})=\varphi(w)$,
hence
$$
   k\,\overline{\varphi(u)c\varphi(u)^{-1}}
     =\overline{\varphi(w)}=p\bar c+\bar w.
$$
Therefore $k=1$ because $p\bar c+\bar w$ cannot be a multiple of another element
of $\K_4^\ab$.
Notice that $\tilde\sigma_1(c)=c$, $\tilde\sigma_1(d^2)=d^2$,
and $\tilde\sigma_1(w)=c^{-1}w$ (see (\eqGLsig)).
Hence, for $f=\tilde\sigma_1^p$, we have
$$
   f\varphi(c)=c,\qquad f\varphi(d^2)=d^2,\qquad\overline{f\varphi(w)}=\bar w.
      \eqno(\eqId)
$$

It remains to show that $Rf\varphi=R\id_{[\tau]}$ for some $\tau:\B'_4\to Z(\B_4)$.
Let $x\in\B'_4$.
Since $\B'_4=\K_4\rtimes\B'_3$ and $\B_4=\K_4\rtimes\B_3$, we may write 
$x=x_1a_1$ and $f\varphi(x)=x_2a_2$ with $x_1=R(x)\in\B'_3$,
$x_2=Rf\varphi(x)\in\B_3$,
and $a_1,a_2\in\K_4$.
The equation (\eqPXP) for $f\varphi$ (and hence with the identity matrix for $P$
because (\eqId) means that $f\varphi|_{\K_4}$ induces the idenity maping of $\K_4^\ab$)
reads $\pi f\varphi(x)=\pi(x)$, that is $\pi(x_2 a_2)=\pi(x_1 a_1)$.
Since $a_1,a_2\in\K_4\subset\ker\pi$, this implies that
$$
     \pi(x_1)=\pi(x_2).                                     \eqno(\eqLemUT)
$$
Let
$S_1=\left(\smallmatrix 1 & -1\\0 & 1\endsmallmatrix\right)$ and
$S_2=\left(\smallmatrix 1 &  0\\1 & 1\endsmallmatrix\right)$.
It is well-known that the mapping $\sigma_1\mapsto S_1$, $\sigma_2\mapsto S_2$
defines an isomorphism between $\B_3/\langle\Delta_3^4\rangle$ and $\SL(2,\Z)$.
From (\eqGLsig) we see that $\pi(\sigma_1)=S_1$ and
$\pi(\sigma_1^{-1}\sigma_2\sigma_1)=S_2$. Hence
$\ker(\pi|_{\B_3})=\langle\Delta_3^4\rangle=R(Z(\B_4))$.
Therefore (\eqLemUT) implies that $x_2 = x_1R(\tau(x))$
for some element $\tau(x)$ of $Z(\B_4)$. It is easy to check that $\tau$
is a group homomorphism, thus, recalling that $x_1=R(x)$ and $x_2=Rf\varphi(x)$,
we get $Rf\varphi(x)=x_2=x_1R(\tau(x))=R(x\tau(x))=R\id_{[\tau]}(x)$.
\qed\enddemo

\proclaim{ Lemma~\lemR } If $\varphi|_{\K_4}=\id$ and
$R\varphi=R\id_{[\tau]}$ for some homomorphism $\tau:\B'_4\to Z(\B_4)$,
then $\varphi=\id_{[\tau]}$.
\endproclaim

\demo{ Proof } Since $\B'_4=\K_4\rtimes\B'_3$ and $\K_4\subset\ker\tau$,
it is enough to show
that $\varphi|_{\B'_3}=\id_{[\tau]}$. So, let $x\in\B'_3$. The condition
$R\varphi=R\id_{[\tau]}$ means that $\varphi(x)=xa\tau(x)$ with $a\in\K_4$.
Let $b$ be any element of $\K_4$. Then $xbx^{-1}\in\K_4$,
hence $\varphi(xbx^{-1})=xbx^{-1}$ (because $\varphi|_{\K_4}=\id$).
Since $\varphi(x)=xa\tau(x)$, $\varphi(b)=b$,
and $\tau(x)$ is central,
it follows that 
$$
  xbx^{-1}=\varphi(xbx^{-1})=\varphi(x)b\varphi(x)^{-1}
      =xa\tau(x)b\tau(x)^{-1}a^{-1}x^{-1} = xaba^{-1}x^{-1}
$$
whence $aba^{-1}=b$. This is true for any $b\in\K_4$, thus
$a\in Z(\K_4)$. Since $\K_4$ is free, we deduce that $a=1$, hence
$\varphi(x)=x\tau(x)=\id_{[\tau]}(x)$.
\qed\enddemo

\demo{ Proof of Theorem \thBBfour }
Recall that we assume in this section that
$\varphi$ is a homomorphism $\B'_4\to\B_4$ such that
$\K_4\not\subset\ker\varphi$.

By Lemma~\lemUT\ we may assume that $\varphi(c)=c$, $\varphi(d^2)=d^2$, and
$R\varphi=R\id_{[\tau]}$ for some $\tau:\B'_4\to Z(\B_4)$,
in particular, $R\varphi(u)=R(u\tau(u))$.
The latter condition means that $\varphi(u)=ua\tau(u)$ with $a\in\K_4$.
Then, by (\eqGL), we have
$$
  \varphi(w)=\varphi(ucu^{-1})=uaca^{-1}u^{-1}=\tilde u(c\,[c^{-1},a]),
      =w[w^{-1},\tilde u(a)],
$$
thus $w^{-1}\varphi(w)=[w^{-1},A]$ for $A=\tilde u(a)\in\K_4$.
By Lemma~\lemWF\ we have also $w^{-1}\varphi(w)\in F$. Then Lemma~\lemXF\
implies that $w^{-1}\varphi(w)=[w^{-1},c^k]$ for some integer $k$, that is
$\varphi(w)=c^k w c^{-k}$. Hence, $(\tilde c^{-k}\varphi)|_{\K_4}=\id$.
Since $c\in\ker R$, we have $R\tilde c^{-k}=R$ whence
$R\tilde c^{-k}\varphi=R\varphi=R\id_{[\tau]}$. This fact combined with
$(\tilde c^{-k}\varphi)|_{\K_4}=\id$ and Lemma~\lemR\ implies that
$\tilde c^{-k}\varphi=\id_{[\tau]}$, i.e., $\varphi$ is equivalent to
$\id_{[\tau]}$.
\qed\enddemo

\Refs
\def\r{\ref}

\r\no\refArtin
\by E.~Artin \paper Theory of braids \jour Ann. of Math. \vol 48 \yr 1947
\pages 101--126
\endref

\r\no\refBM
\by R.~W.~Bell, D.~Margalit
\paper Braid groups and the co-Hopfian property
\jour J. Algebra \vol 303 \yr 2006 \pages 275--294
\endref


\r\no\refBLM
\by J.~S.~Birman, A.~Lubotzky, J.~McCarthy
\paper Abelian and solvable subgroups of the mapping class group
\jour Duke Math. J. \vol  50 \yr 1983 \pages 1107--1120
\endref


\ref\no\refCastel
\by F.~Castel
\paper Geometric representations of the braid groups
\jour Ast\'erisque \vol 378 \yr 2016 \pages vi+175
\endref


\r\no\refDG
\by    J.~L.~Dyer, E.~K.~Grossman
\paper The automorphism group of the braid groups
\jour  Amer. J. of Math. \vol 103 \yr 1981 \pages 1151--1169
\endref

\r\no\refFM
\by B.~Farb, D.~Margalit
\book A primer on mapping class groups
\bookinfo volume 49 of Princeton Mathematical Series
\publ Princeton University Press \publaddr Princeton, NJ \yr 2012
\endref

\r\no\refGM
\by    J.~Gonz\'alez-Meneses
\paper The $n$th root of a braid is unique up conjugacy
\jour  Algebraic and Geometric Topology \vol 3 \yr 2003 \pages 1103--1118
\endref

\r\no\refGMW
\by    J.~Gonz\'alez-Meneses, B.~Wiest
\paper On the structure of the centralizer of a braid
\jour  Ann. Sci. \'Ec. Norm. Sup\'er. (4) \vol 37 \yr 2004 \pages 729--757
\endref

\r\no\refGL
\by    E.~A.~Gorin, V.~Ya.~Lin
\paper Algebraic equations with continuous coefficients and some problems
       of the algebraic theory of braids
\jour  Math. USSR-Sbornik \vol 7 \yr 1969 \pages 569--596.
\endref

\r\no\refI
\by    N.~V.~Ivanov
\book  Subgroups of Teichm\"uller modular groups
\bookinfo Translations of mathematical monographs \vol 115 \yr 1992 \publ AMS
\endref

\r\no\refKM
\by    K.~Kordek, D.~Margalit
\paper Homomorphisms of commutator subgroups of braid groups
\jour 	arXiv:1910.06941
\endref

\r\no\refLin
\by    V.~Lin \paper Braids and permutations \jour arXiv:math/0404528 \endref

\ref\no\refLinFAA
\by    V.~Ya.~Lin
\paper Algebraic functions, configuration spaces,
       Teichm\"uller spaces, and new holomorphically combinatorial invariants
\jour Funk. Anal. Prilozh. \vol 45 \yr 2011 \issue 3 \pages 55--78
\lang Russian \transl English transl.
\jour Funct. Anal. Appl. \vol 45 \yr 2011 \issue 3 \pages 204--224
\endref

\r\no\refLinTalk
\by    V.~Lin
\paper Some problems that I would like to see solved
\jour  Abstract of a talk. Technion, 2015,
http://www2.math.technion.ac.il/$\widetilde{\;}$pincho/Lin/Abstracts.pdf
\endref

\r\no\refMKS
\by    W.~Magnus, A.~Karrass, D.~Solitar
\book  Combinatorial group theory: presentations of groups in terms of generators and relations
\publ  Interscience Publ. \yr 1966
\endref

\r\no\refJAlg
\by     S.~Yu.~Orevkov
\paper  Algorithmic recognition of quasipositive braids of algebraic length two
\jour   J. of Algebra \vol 423 \yr 2015 \pages 1080--1108
\endref

\r\no\refAFST
\by     S.~Yu.~Orevkov
\paper  Automorphism group of the commutator subgroup of the braid group
\jour   Ann. Facult\'e des Scie. de Toulouse. Math. (6)
        \vol 26 \yr 2017 \pages 1137--1161
\endref


\endRefs
\enddocument